\documentclass[12pt]{amsart}
\usepackage{latexsym, amsmath, amssymb, amsthm}
\usepackage{epsfig}
\usepackage{amscd}
\usepackage{latexsym}
\usepackage{pstricks,pst-node,cite}

\newtheorem{thm}{Theorem}

\newcommand{\RR}{\ensuremath{\mathbb{R}^2}}

\newcommand{\Inv}{\mbox{\rm Inv}}

\newcommand{\uu}{{\uparrow}}
\newcommand{\dd}{{\downarrow}}

\newcounter{fignum}
\ifx\blackandwhite\undefined
  \psset{linecolor=blue}\newrgbcolor{darkred}{.9 0 0}\newrgbcolor{emgreen}{0 .9 0}
  \newrgbcolor{pup}{.7  0 .9}\newrgbcolor{xxx}{.8 .8 .8} 
    \newrgbcolor{lightgray}{0 .9 0}\newrgbcolor{darkgray}{.99 .99 0} \else
  \psset{linecolor=black}\newgray{darkred}{.7}\newgray{emgreen}{0.3}\newgray{pup}{.5}\newgray{xxx}{.5}
\fi

\psset{linewidth=.4pt,dash=3pt 3pt,doublesep=.05,dotsize=1pt 5}
\SpecialCoor

\begin{document}

\author{Dongseok Kim}
\address{Department of Mathematics \\Kyonggi University
\\ Suwon, 443-760 Korea}
\email{dongseok@kgu.ac.kr}
\thanks{}



\title[On the $(n,k)$-th Catalan numbers]
{On the $(n,k)$-th Catalan numbers}
\begin{abstract}
In this paper, we generalize the Catalan number to the $(n,k)$-th Catalan numbers and find
a combinatorial description that the $(n,k)$-th Catalan numbers is equal to the number of partitions of
$n(k-1)+2$ polygon by $(k+1)$-gon where all vertices of all $(k+1)$-gons lie on the vertices of
$n(k-1)+2$ polygon.
\end{abstract}

\maketitle

\section{Introduction} \label{intro}

Let $D^2$ be the closed unit disc in $\RR$. We label $2n$ points on
the boundary of $D^2$ by $\{ 1$, $2$, $\ldots$, $2n\}$ with a counterclockwise
orientation. We embed $n$ arcs into $D^2$ so that the boundaries
of intervals map to the labeled points. Such an embedding is called an
\emph{$(n,2)$ diagram} if it does not have any crossings. We consider the isotopy
classes of $(n,2)$ diagrams relative to the boundary. The number of isotopy
classes of $(n,2)$ diagrams is called the \emph{$n$-th Catalan number}~\cite{Catalan}, denoted by $C_n$ and
$$C_n = \frac{1}{n+1} \binom {2n}{n}.$$
$66$ different combinatorial descriptions of
$n$-th Catalan number are given in a famous literature~\cite[Exercise 6. 19]{Stanley:em2}
and currently $96$ more are added on his Catalan addendum~\cite{Stanley:home}. This number
plays a major rule in many different areas of not only mathematics but also natural science~\cite{CG, Gardner}.
Harary found an implication of $n$-th Catalan number on the cell growth problem~\cite{Harary}.
Cameron \cite{Cameron} found representations for binary and ternary trees which can be used to enumerate Rothe numbers~\cite{Gould:binomial}.

\begin{figure}
$$
\begin{pspicture}[shift=-2.4](-2.5,-2.5)(2.5,2.5)
\pnode(2;0){a1} \pnode(2;45){a2} \pnode(2;90){a3} \pnode(2;135){a4}
\pnode(2;180){a5} \pnode(2;225){a6}\pnode(2;270){a7} \pnode(2;315){a8}
\rput(2.3;0){$1$} \rput(2.3;45){$2$} \rput(2.3;90){$3$} \rput(2.3;135){$4$}
\rput(2.3;180){$5$} \rput(2.3;225){$6$} \rput(2.3;270){$7$} \rput(2.3;315){$8$}
\pscircle[linewidth=1pt, linestyle=dashed](0,0){2}
\nccurve[angleA=-90,angleB=-45]{a3}{a4}
\nccurve[angleA=0,angleB=-135]{a5}{a2}
\nccurve[angleA=45,angleB=180]{a6}{a1}
\nccurve[angleA=90,angleB=135]{a7}{a8}
\end{pspicture} \quad\quad\quad\quad\quad \begin{pspicture}[shift=-2.4](-2.5,-2.5)(2.5,2.5)
\pnode(2;0){a1} \pnode(2;20){a2} \pnode(2;40){a3} \pnode(2;60){a4}
\pnode(2;80){a5} \pnode(2;100){a6}\pnode(2;120){a7} \pnode(2;140){a8}
\pnode(2;160){a9} \pnode(2;180){a10} \pnode(2;200){a11} \pnode(2;220){a12}
\pnode(2;240){a13} \pnode(2;260){a14}\pnode(2;280){a15} \pnode(2;300){a16}
\pnode(2;320){a17} \pnode(2;340){a18}
\pnode(1;0){b1} \pnode(1;20){b2} \pnode(1;40){b3} \pnode(1;60){b4}
\pnode(1;80){b5} \pnode(1;100){b6}\pnode(1;120){b7} \pnode(1;140){b8}
\pnode(1;160){b9} \pnode(1;180){b10} \pnode(1;200){b11} \pnode(1;220){b12}
\pnode(1;240){b13} \pnode(1;260){b14}\pnode(1;280){b15} \pnode(1;300){b16}
\pnode(1;320){b17} \pnode(1;340){b18} \pnode(0;0){c0} \pnode(.5;130){d0}
\rput(2.3;0){$1$} \rput(2.3;20){$2$} \rput(2.3;40){$3$} \rput(2.3;60){$4$}
\rput(2.3;80){$5$} \rput(2.3;100){$6$} \rput(2.3;120){$7$} \rput(2.3;140){$8$}
\rput(2.3;160){$9$} \rput(2.3;180){$10$} \rput(2.3;200){$11$} \rput(2.3;220){$12$}
\rput(2.3;240){$13$} \rput(2.3;260){$14$} \rput(2.3;280){$15$} \rput(2.3;300){$16$}
\rput(2.3;320){$17$} \rput(2.3;340){$18$}
\pscircle[linewidth=1pt, linestyle=dashed](0,0){2}  \ncline{c0}{a17}
\ncline{b1}{a1} \ncline{b5}{a5} \ncline{b10}{a10} \ncline{b15}{a15}
\nccurve[angleA=50,angleB=-140]{c0}{a3}
\nccurve[angleA=-130,angleB=60]{c0}{a13}
\nccurve[angleA=-160,angleB=90]{a2}{b1}
\nccurve[angleA=-90,angleB=160]{b1}{a18}
\nccurve[angleA=-80,angleB=170]{a6}{b5}
\nccurve[angleA=-10,angleB=-120]{b5}{a4}
\nccurve[angleA=-20,angleB=90]{a9}{b10}
\nccurve[angleA=-90,angleB=20]{b10}{a11}
\nccurve[angleA=80,angleB=-170]{a14}{b15}
\nccurve[angleA=10,angleB=120]{b15}{a16}
\nccurve[angleA=-60,angleB=40]{a7}{d0}
\nccurve[angleA=-40,angleB=-140]{a8}{d0}
\nccurve[angleA=40,angleB=-50]{a12}{d0}
\end{pspicture}
$$
\caption{A $(4,2)$ diagram and a $(6,3)$ diagram.}\label{figexample1}
\end{figure}
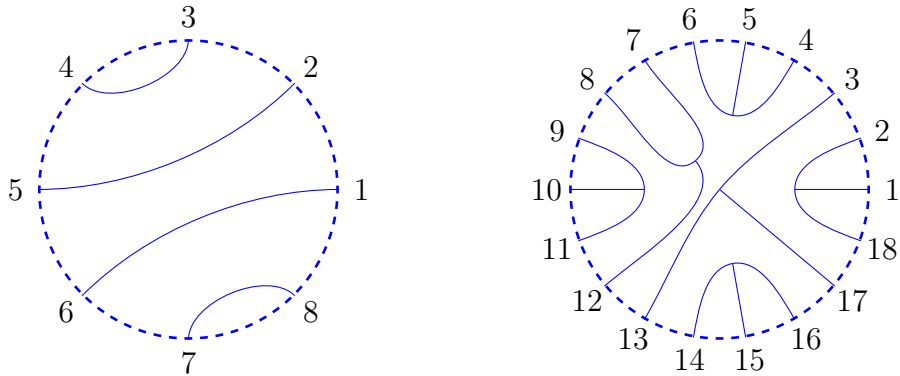

There have been many attempts to generalize Catalan numbers. We are particularly interested in two of them.
Gould developed a generalization of $n$-th Catalan number as follows~\cite{Gould:identity},
$$A_n(a,b) = \frac{a}{a+bn} \binom {a+bn}{n},$$
together with its convolution formula,
$$\sum_{k=0}^{n} A_k(a,b)A_{n-k}(a,b) = A_n(a+c,b).$$
From Gould's generalization, $A_n(a,b)$, one can obtain $C_n$ by substituting $a=1$ and $b=2$.
On the otherhand, $n$-th Catalan number is the dimension of the linear skein space of a disc with $2n$
points on its boundary or simply $n$-th Temperley-Lieb algebras, $\mathcal{T}_n$. Kuperberg generalized
Temperley-Lieb algebras, which corresponds to the invariant subspace
of a tensor product of the vector representation of
${\mathfrak{sl}}(2)$, to web spaces of simple Lie algebras of rank
$2$, ${\mathfrak{sl}}(3)$, ${\mathfrak{sp}}(4)$ and $G_2$
\cite{Kuperberg:spiders}. It has been extensively studied~\cite{Kim:tri, Kim:JW, KL:sl3}.
In particular, a purpose of the present article is to find a generalization which fits in both directions;
Gould's generalization and the web spaces of ${\mathfrak{sl}}(3)$.
The following describes our generalization. For
$k\ge 2, n\ge 1$, we label $kn$ points on the boundary of the disc by $\{ 1, 2, \ldots, kn\}$.
Then we embed n copies of $k$-stars with  boundaries of $k$-stars
goes the labeled points where a $k$-star is a neighborhood of a
vertex in the complete graph $K_{k+1}$ on $k+1$ vertices. Then we consider  an $(n,k)$
diagram and its isotopy class in a similar fashion. A $(4,2)$ diagram and a $(6,3)$ diagram are illustrated in Figure~\ref{figexample1}.
Let $\mathcal{D}(n,k)$ be the set of all isotopy classes of $(n,k)$-diagrams. The number of elements in  $\mathcal{D}(n,k)$
is called the \emph{$(n, k)$-Catalan number}.

\begin{thm}
The $(n, k)$-th Catalan number is $\frac{1}{n} \binom {kn}{n-1}$. \label{thm1}
\end{thm}

The first combinatorial description of
$n$-th Catalan number given in~\cite[Exercise 6. 19~{\bf a}]{Stanley:em2} is that it is equal to
the number of triangulations of $(n+2)$ polygon $P_{n+2}$ where all vertices of triangles lie on the
vertices of $P_{n+2}$. We also find a similar description for $(n, k)$-Catalan number
in the following theorem.

\begin{thm}
The $(n,k)$-th Catalan number is equal to the number of partitions
of $n(k-1)+2$ polygon by $(k+1)$-gon where all vertices of all
$(k+1)$-gon lie on the vertices of $n(k-1)+2$ polygon. \label{thm2}
\end{thm}

The outline of this paper is as follows. In section~\ref{proof}, we
provide the proof of Theorem~\ref{thm1} and~\ref{thm2}. In section~\ref{diss} we
illustrate how $(n, k)$-th Catalan number is related with other generalizations and
the representation theory of the quantum Lie algebra $\mathcal{U}_q(\mathfrak{sl}(3,\mathbb{C}))$.

\section{Proof of main theorems} \label{proof}

\subsection{Proof of Theorem~\ref{thm1}}

Let $I_{l}=\{ 1, 2, \cdots, l\}$ be $\mathbb{Z}/l \mathbb{Z}$.
Let $\mathcal{S}(l,m)$ be the set of all subsets of $I_l$ whose
cardinalities are $m$. We first observe the
relation between $\mathcal{S}(nk, n)$ and  $\mathcal{D}(n,k)$ by defining
a map $\psi : \mathcal{S}(nk, n) \rightarrow \mathcal{D}(n,k)$. Let
$E$ be an element in $\mathcal{S}(nk,n)$. To obtain the image $\psi (E)$, an
$(n,k)$ diagram, we first array $nk$ points on the boundary of the unit
disc $D^2$. We assign $\uu$ for $n$ points in $E$, $\dd$ for
points in $I_{nk}-E$. Let $j$ be the smallest number in $E$. Now we
define a height function $\varphi : I_{nk} \rightarrow{\mathbb
Z}$ by $\varphi (i) =$ the number of $\uu$'s minus the number of
$\dd$'s  from $j$ to $i$ in the clockwise orientation on the
boundary of the disc $D^2$. Then clearly $\varphi (j)=1$, $\varphi
(j-1)=-(k-1)n$ and $\varphi$ is a step function such that $\vert
\varphi(l)-\varphi(l+1)\vert =1$ for all $l=1, 2, \cdots, kn$ except
$j-1$. Since it starts from $1$ at $j$ and ends with $-(k-1)n$ at
$j-1$, by a simple application of the Pigeon Hole principle, we can see that there exists $l$ such that

$$\varphi (l-1) < \varphi (l)  = \varphi (l+1)+1 =\cdots =\varphi (l+k-1) +k-1,$$
which give us an innermost $k$-star joining $k$ consecutive points
from $l$. After removing these $k$ consecutive points, inductively
we repeat this process from the next innermost $k$-star. By picking
the first point for the points of each $k$-star from $j$, we can see
that $\psi$ is onto. For example, we will pick $\{ 3, 4, 7, 9, 14,
18\}$ for $(6,3)$ diagram in the righthand side of Figure~\ref{figexample1} and we will get this
$(6,3)$ diagram from this set of six points if we follow the inductive
process.

One can easily notice that there are different sets of $n$ points
which make the same $(n,k)$ diagram, $i.e.$, $\psi$ is not one to one, such
as $\{4, 7, 9, 13, 14, 18\}$ and $\{ 1, 3, 4, 7, 9, 14\}$ will give the same $(6,3)$ diagram
in Figure~\ref{figexample1}. However, it is not clear to
find how many different subsets in $\mathcal{S}(nk, n)$ lead us to the same $(n,k)$
diagram. To get the right number of all distinct $(n,k)$ diagrams, we will
show the second claim : there exists a surjective map $\vartheta : \mathcal{
S}(nk, n-1) \rightarrow \mathcal{D}(n,k)$, $i.e.$, for a given element $F$ of $\mathcal{
S}(nk,n-1)$, there is a unique way to decide the last point we have to pick so that we get an element
$E$ in $\mathcal{S}(nk,n)$ and an $(n,k)$ diagram, $\varphi (E) =
\vartheta (F)$. The proof follows from the proof of the first claim
because one can see that at the last stage of the induction procedure, we had $k$ points
left and the $(n,k)$ diagram was decided already. Therefore, we just pick the
smallest label from these $k$ points.
For the surjectivity of $\vartheta$, for a given $(n,k)$ diagram $D$,
we knew that there is an innermost $k$-star from
the first claim. Then to obtain the desired set we first start to read the first label for each $k$
star from the last label of the innermost $k$-star. Then we
will have an ordered set $E$ of $n$ points and furthermore, if we drop the last
label to get an ordered set $F$ of $n-1$ points, the previous process will recover
it, $i. e.$, $\varphi (E) = \vartheta (F)=D$.

The last claim for the proof is that $\vartheta$ is an $n$ to $1$ map, $i. e.$, for a
given an $(n,k)$ diagram there are $n$ unique sets of $n-1$ points in $\mathcal{
S}(nk, n-1)$ such that if we apply the process in the second claim, we recover
the original diagram. For a given $(n,k)$ diagram $D$ and each $k$ star in $D$, we drop
it and read the first label in other $(n-1)$ $k$-stars. One can
easily see that all these $(n-1)$ points will recover the same $(n,k)$ diagram
by the process described in the second claim. Also we can see that if we pick
a point other than these $(n-1)$ points, then we may not recover the
$k$ star we have dropped. For example, $\{3, 4, 7, 9, 14\}$, $\{ 4, 7,
9, 14, 18\}$, $\{4$, $7$, $9$, $ 7$, $18\}$, $\{4, 9, 13, 14, 18\}$, $\{4,
12, 13, 14, 18\}$ and $\{7, 9, 13 ,14, 18\}$ are only possible six sets
of five elements in $\mathcal{S}(18, 5)$ for the $(6,3)$ diagram in
Figure~\ref{figexample1}. Therefore, the $(n,k)$ Catalan number
is the one $n$-th of the cardinality of $\mathcal{S}(nk, n-1)$ which is $\frac{1}{n}\binom{kn}{n-1}$.

\subsection{Proof of Theorem~\ref{thm2}}

\begin{figure}
$$
\begin{pspicture}[shift=-.4](.8,-.5)(9.8,3.8)
\psline[linewidth=1pt, linestyle=dashed](.7,0)(9.8,0)
\rput(9,-.3){$1$} \rput(9.5,-.3){$2$} \rput(1,-.3){$3$}
\rput(1.5,-.3){$4$} \rput(2,-.3){$5$} \rput(2.5,-.3){$6$}
\rput(3,-.3){$7$} \rput(3.5,-.3){$8$} \rput(4,-.3){$9$}
\rput(4.5,-.3){$10$} \rput(5,-.3){$11$} \rput(5.5,-.3){$12$}
\rput(6,-.3){$13$} \rput(6.5,-.3){$14$} \rput(7,-.3){$15$}
\rput(7.5,-.3){$16$} \rput(8,-.3){$17$} \rput(8.5,-.3){$18$}
\psarc(2,0){.5}{0}{180} \psline(2,0)(2,.5)
\psarc(4.5,0){.5}{0}{180} \psline(4.5,0)(4.5,.5)
\psarc(7,0){.5}{0}{180} \psline(7,0)(7,.5)
\psarc(4.25,0){1.25}{0}{180} \psarc(4.75,0){1.25}{101.5}{180}
\psarc(4.5,0){3.5}{0}{180} \psarc(2.5,0){3.5}{0}{73.5}
\psarc(9,0){.5}{0}{180} \psline(9,0)(9,.5)
\end{pspicture}
$$
\caption{The $(6,3)$ diagram $\tilde D$ in the upperhalf plane corresponding to
a $(6,3)$ diagram $D$ in Figure~\ref{figexample1}.}\label{figexample2}
\end{figure}
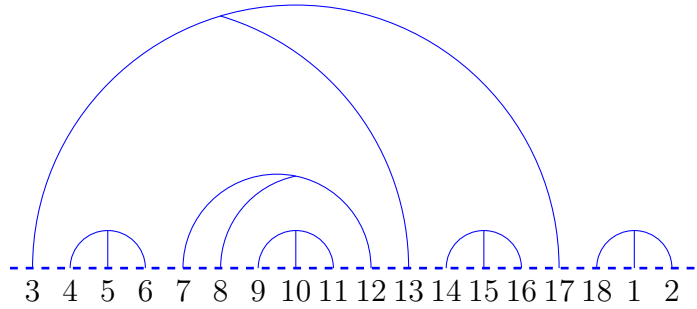

The main idea of the proof came from that of the proof for the Catalan number by
replacing $\uu$ for the first left branch and $\dd$ for the other
branches~\cite{Stanley:em2}.

\begin{figure}
$$
\begin{pspicture}[shift=-.4](-.2,-.5)(12.2,4.5)
\rput(6,4.4){$T_D$}
\rput(1.8,3.2){$a$} \rput(8.2,3.2){$b$} \rput(11.2,3.2){$c$}
\rput(0,1.7){$aa$} \rput(1,1.7){$ab$} \rput(4.3,2.2){$ac$}
\rput(7,1.7){$ba$} \rput(8,1.7){$bb$} \rput(9,1.7){$bc$}
\rput(10,1.7){$ca$} \rput(11,1.7){$cb$} \rput(12,1.7){$cc$}
\rput(2,.7){$aca$} \rput(3.6,1.3){$acb$} \rput(6,.7){$acc$}
\rput(3,-.3){$acba$} \rput(4,-.3){$acbb$} \rput(5,-.3){$acbc$}
\psline(6,4)(2,3)\psline(6,4)(8,3)\psline(6,4)(11,3)
\psline(2,3)(0,2)\psline(2,3)(1,2)\psline(2,3)(4,2)
\psline(8,3)(7,2)\psline(8,3)(8,2)\psline(8,3)(9,2)
\psline(11,3)(10,2)\psline(11,3)(11,2)\psline(11,3)(12,2)
\psline(4,2)(2,1)\psline(4,2)(4,1)\psline(4,2)(6,1)
\psline(4,1)(3,0)\psline(4,1)(4,0)\psline(4,1)(5,0)
\pscircle[fillstyle=solid,fillcolor=darkgray,linecolor=black](6,4){.13}
\pscircle[fillstyle=solid,fillcolor=darkgray,linecolor=black](2,3){.13}
\pscircle[fillstyle=solid,fillcolor=darkgray,linecolor=black](8,3){.13}
\pscircle[fillstyle=solid,fillcolor=darkgray,linecolor=black](11,3){.13}
\pscircle[fillstyle=solid,fillcolor=darkgray,linecolor=black](0,2){.13}
\pscircle[fillstyle=solid,fillcolor=darkgray,linecolor=black](1,2){.13}
\pscircle[fillstyle=solid,fillcolor=darkgray,linecolor=black](4,2){.13}
\pscircle[fillstyle=solid,fillcolor=darkgray,linecolor=black](7,2){.13}
\pscircle[fillstyle=solid,fillcolor=darkgray,linecolor=black](8,2){.13}
\pscircle[fillstyle=solid,fillcolor=darkgray,linecolor=black](9,2){.13}
\pscircle[fillstyle=solid,fillcolor=darkgray,linecolor=black](10,2){.13}
\pscircle[fillstyle=solid,fillcolor=darkgray,linecolor=black](11,2){.13}
\pscircle[fillstyle=solid,fillcolor=darkgray,linecolor=black](12,2){.13}
\pscircle[fillstyle=solid,fillcolor=darkgray,linecolor=black](2,1){.13}
\pscircle[fillstyle=solid,fillcolor=darkgray,linecolor=black](4,1){.13}
\pscircle[fillstyle=solid,fillcolor=darkgray,linecolor=black](6,1){.13}
\pscircle[fillstyle=solid,fillcolor=darkgray,linecolor=black](3,0){.13}
\pscircle[fillstyle=solid,fillcolor=darkgray,linecolor=black](5,0){.13}
\pscircle[fillstyle=solid,fillcolor=darkgray,linecolor=black](4,0){.13}
\end{pspicture}
$$
\caption{The rooted tree $T_D$ corresponding to the $(6,3)$ diagram $D$ in Figure~\ref{figexample1}
with words describe the lexicographic order to recover the $(6,3)$ diagram $D$.}\label{figexample3}
\end{figure}
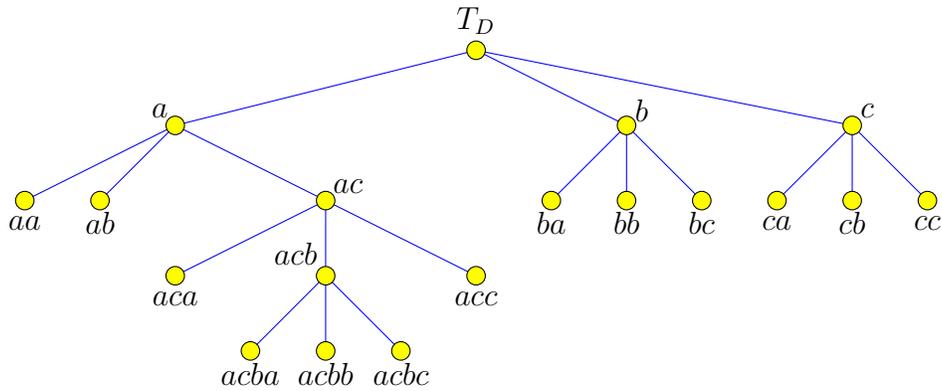

Let $D$ be an $(n,k)$ diagram. We pick a set $E$ of $n$ points as we
described in the proof for the surjectivity of $\varphi$. We cut open the
disc $D^2$ from a point $p$ which is in between the smallest label in $E$
and its immediate predecessor so that we get a diagram $\tilde D$ on the upperhalf plane
as in Figure~\ref{figexample2}.
One can notice that all points in $E$, for $\tilde D$ in Figure~\ref{figexample2} $E$ can be chosen
$\{3, 4, 7, 9, 14, 18\}$, appears in the leftmost leg of
each $k$-star. Then we can make a rooted
$k$-ary tree $T_D$ by the following way. First we start from the $k$-star
contains the smallest label in $E$. We just draw $k$ branches from the root. Since each
$k$-star divides the upper half plane into $k$ disjoint regions and
we have the counterclockwise orientation around the \emph{center} of
$k$ star, where all legs of $k$-star join, it gives us a cyclic
order of the edges of $k$ star. So we can order the regions from the next region of the unbounded region
in counterclockwise orientation. If a region does not
contain any other $k$-stars, that branch stops. Otherwise, we repeat the
process from the region in the leftmost $k$-star, at
this stage we ignore the previous $k$-stars. Inductively, we can
complete the rooted tree $\tilde T$. For converse, from the root
we assign alphabet letters for children nods from the left.
Inductively at each generation of parents nods, for each children nods we append
alphabet letters to the word assigned to the parent nod from the left. So, it assigns
words for every nods except the root. Then we can array them in line by the lexicographic
dictionary order. For the example in Figure~\ref{figexample3}, we obtain the following words for each labels:
$3(a), 4(aa), 5(ab), 6(ac)$, $7(aca)$, $8(acb)$, $9(acba)$, $10(acbb)$, $11(acbc)$, $12(acc)$, $13(b)$, $14(ba)$,
$15(bb)$, $16(bc)$, $17(c)$, $18(ca)$, $1(cb)$, $2(cc)$. Then we read $n$ points which is the leftmost branch
or words which end with the letter $a$ to recover $n$ points set $E$. One can see this process in Figure~\ref{figexample3}.

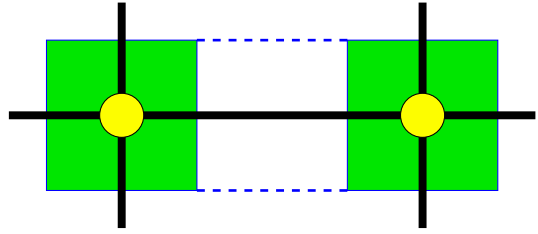
\begin{figure}
$$
\begin{pspicture}[shift=-.4](-1.7,-1.7)(5.7,1.7)
\psline[linewidth=1pt, linestyle=dashed](1,1)(3,1)
\psline[linewidth=1pt, linestyle=dashed](1,-1)(3,-1)
\pspolygon[fillstyle=solid, fillcolor=lightgray](1,1)(-1,1)(-1,-1)(1,-1)(1,1)
\pspolygon[fillstyle=solid, fillcolor=lightgray](5,1)(3,1)(3,-1)(5,-1)(5,1)
\psline[linecolor=black, linewidth=3pt](-1.5,0)(5.5,0)
\psline[linecolor=black, linewidth=3pt](0,-1.5)(0,1.5)
\psline[linecolor=black, linewidth=3pt](4,-1.5)(4,1.5)
\pscircle[fillstyle=solid,fillcolor=darkgray,linecolor=black](0,0){.3}
\pscircle[fillstyle=solid,fillcolor=darkgray,linecolor=black](4,0){.3}
\end{pspicture}
$$
\caption{Strips for connecting the adjacent quadrilaterals from the rooted tree $T_D$ and dotted lines which
will be collapsing to a point to get $n(k-1)+2$ polygon.}\label{figexample4}
\end{figure}

\begin{figure}
$$ \begin{pspicture}[shift=-2.4](-4,-4)(4,4)
\pnode(3;0){a1} \pnode(3;25.7){a2} \pnode(3;51.4){a3}
\pnode(3;77){a4} \pnode(3;103){a5} \pnode(3;128.6){a6} \pnode(3;154.3){a7} \pnode(3;180){a8} \pnode(3;205.5){a9} \pnode(3;231){a10} \pnode(3;256.5){a11} \pnode(3;282){a12}
\pnode(3;308){a13} \pnode(3;334){a14}\pnode(3;360){a15}
\psline(3.4;347)(1;38.5)
\psline(3.4;12.8)(1;38.5)(1.6;90)(2.6;90)(3.4;90)
\psline(3.4;243.8)(2.6;270)(3.4;295)
\psline(3.4;166.5)(2.6;192.7)(3.4;218.3)
\psline(3.4;38.5)(1.6;90)(3.4;141.4)
\psline(3.4;270)(2.6;270)(1;270)(1;38.5)
\psline(3.4;64.2)(2.6;90)(3.4;115.3)
\psline(1;270)(2.6;192.7)(3.4;192.7)
\ncline[linewidth=1pt]{a1}{a2}
\ncline[linewidth=1pt]{a3}{a2}
\ncline[linewidth=1pt]{a3}{a4}
\ncline[linewidth=1pt]{a5}{a4}
\ncline[linewidth=1pt]{a5}{a6}
\ncline[linewidth=1pt]{a7}{a6}
\ncline[linewidth=1pt]{a7}{a8}
\ncline[linewidth=1pt]{a9}{a8}
\ncline[linewidth=1pt]{a9}{a10}
\ncline[linewidth=1pt]{a11}{a10}
\ncline[linewidth=1pt]{a11}{a12}
\ncline[linewidth=1pt]{a13}{a12}
\ncline[linewidth=3pt]{a13}{a14}
\ncline[linewidth=1pt]{a14}{a15}
\ncline[linewidth=1pt, linestyle=dashed, linecolor=darkred]{a14}{a7}
\ncline[linewidth=1pt, linestyle=dashed, linecolor=darkred]{a2}{a7}
\ncline[linewidth=1pt, linestyle=dashed, linecolor=darkred]{a7}{a10}
\ncline[linewidth=1pt, linestyle=dashed, linecolor=darkred]{a10}{a13}
\ncline[linewidth=1pt, linestyle=dashed, linecolor=darkred]{a3}{a6}
\pscircle[fillstyle=solid,fillcolor=darkgray,linecolor=black](3.4;12.8){.13}
\pscircle[fillstyle=solid,fillcolor=darkgray,linecolor=black](3.4;38.5){.13}
\pscircle[fillstyle=solid,fillcolor=darkgray,linecolor=black](3.4;64.2){.13}
\pscircle[fillstyle=solid,fillcolor=darkgray,linecolor=black](3.4;90){.13}
\pscircle[fillstyle=solid,fillcolor=darkgray,linecolor=black](3.4;115.3){.13}
\pscircle[fillstyle=solid,fillcolor=darkgray,linecolor=black](3.4;141.4){.13}
\pscircle[fillstyle=solid,fillcolor=darkgray,linecolor=black](3.4;166.5){.13}
\pscircle[fillstyle=solid,fillcolor=darkgray,linecolor=black](3.4;192.7){.13}
\pscircle[fillstyle=solid,fillcolor=darkgray,linecolor=black](3.4;218.3){.13}
\pscircle[fillstyle=solid,fillcolor=darkgray,linecolor=black](3.4;243.8){.13}
\pscircle[fillstyle=solid,fillcolor=darkgray,linecolor=black](3.4;270){.13}
\pscircle[fillstyle=solid,fillcolor=darkgray,linecolor=black](3.4;295){.13}
\pscircle[fillstyle=solid,fillcolor=darkgray,linecolor=black](3.4;347){.13}
\pscircle[fillstyle=solid,fillcolor=darkgray,linecolor=black](1.6;90){.13}
\pscircle[fillstyle=solid,fillcolor=darkgray,linecolor=black](2.6;192.7){.13}
\pscircle[fillstyle=solid,fillcolor=darkgray,linecolor=black](2.6;270){.13}
\pscircle[fillstyle=solid,fillcolor=darkgray,linecolor=black](1;38.5){.13}
\pscircle[fillstyle=solid,fillcolor=darkgray,linecolor=black](2.6;90){.13}
\pscircle[fillstyle=solid,fillcolor=darkgray,linecolor=black](1;270){.13}
\rput(3.2;322){$e$}
\rput(3.8;12.8){$ab$} \rput(3.8;38.5){$aca$} \rput(3.8;64.2){$acba$}
\rput(3.8;90){$acbb$} \rput(3.8;115.3){$acbc$} \rput(3.8;141.4){$acc$}
\rput(3.8;166.5){$ba$} \rput(3.8;192.7){$bb$} \rput(3.8;218.3){$bc$} \rput(3.8;243.8){$ca$}
\rput(3.8;270){$cb$} \rput(3.8;295){$cc$} \rput(3.8;347){$aa$}
\end{pspicture}
$$
\caption{The partition of $14$-gon by rectangles obtained from the rooted tree $T_D$ in Figure~\ref{figexample4}.}\label{figexample5}
\end{figure}
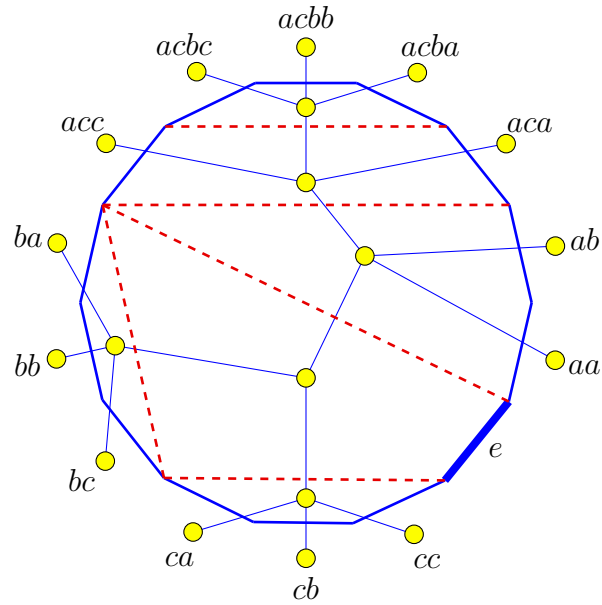

For the next step, we cover all vertices of valence bigger than one by small
quadrilaterals such that they only intersect the tree at the neighborhood
of the center vertex and every edges of quadrilaterals centered at the
vertex are transversal to only one edge of the neighborhood of the
center vertex as depicted in Figure~\ref{figexample4}. We must have exactly $n$ vertices of valence bigger than one.
We make an exception for the root where we allow to have an edge of
quadrilateral centered at root does not intersect any part of the tree.
Then each edges between two vertices of valence bigger than one, we put a
strip which connects the edges of quadrilaterals facing each other as in
Figure~\ref{figexample2}. At this stage, we have $n(k+1)$ polygon. Then we collapse
these extra strips by shrinking each dotted line to a point as depicted in
Figure~\ref{figexample5}. This process leads us to an $n(k-1)+2$ polygon. Finally, we transform it to a regular
$n(k-1)+2$-gon where an edge of the quadrilateral contains the root
which does not intersect the tree $T$ goes to a fixed edge $e$.
The converse part should be clear. It completes the proof.

\section{Discussion and problems} \label{diss}

In this section, we discuss the relations between our generalization of the Catalan number
and other generalizations and the representation theory of quantum Lie algebra
$\mathcal{U}_q(\mathfrak{sl}(3,\mathbb{C}))$.

The first easy observation
is to see that the $(n,k)$-th Catalan number can be obtained from $A_n(a,b)$
by substituting $a=1$ and $b=k$. Consequently, the $n$-th Catalan number
$C_n$ can be obtained from the $(n,k)$-th Catalan number by substituting $k=2$.
Although, we have not stated separate theorems but the proof of Theorem~\ref{thm2}
showed two more combinatorial descriptions of $(n,k)$-th Catalan number by the rooted $k$-ary tree and
alphabetical words~\cite[Exercise 6. 19 d]{Stanley:em2}.

As we stated before, the $n$-th Catalan number is the dimension of the invariant space of $V_1^{\otimes 2n}$ where
$V_1$ is the vector representation of $\mathfrak{sl}(2,\mathbb{C})$ because $(n,2)$ diagrams are precisely
geometric realizations of the invariant vectors. Thus, it is natural to ask $(n,3)$ diagrams may be also
geometric realizations of the invariant vectors in the invariant space of $V_1^{\otimes 3n}$ where
$V_1$ is the vector representation of $\mathfrak{sl}(3,\mathbb{C})$. In fact, it is known that there exists a general method
to generate all invariant vectors in the invariant space of $V_1^{\otimes 3n}$~\cite{KL:sl3} and our expectation is right.
However, not all invariant vectors are $(n,3)$ diagrams. Because the full invariant space can be described as
\begin{align*}
\Inv(V_{1}^{\otimes 3n})
&= \bigoplus_{k=1}^{3n-2}\bigoplus_{V_{w_1}\in \mathcal{C}} \ldots \bigoplus_{V_{w_{k}}\in
\mathcal{C}} \Inv(V_{1}^{\otimes a_1} \otimes V_{w_1}))\oplus
\Inv(V_{w_1}^* \otimes V_{1}^{\otimes a_2} \otimes V_{w_2})\\ &\oplus\ldots
\oplus\Inv(V_{w_{k}}^* \otimes V_{1}^{\otimes a_k}).\label{tensordecomp}
\end{align*}
where $\mathcal{C}$ be the set of all irreducible representations of
${\mathfrak{sl}}(3,\mathbb{C})$. But the invariant vectors corresponding to $(n,3)$ diagrams are the
subspace of $\Inv(V_{1}^{\otimes 3n})$ where $V_{w_i}$ is the trivial representation, $\mathbb{C}$.

\vskip 1cm
\noindent{\bf Acknowledgements}
The author would like to thank Greg Kuperberg for introducing the
subject, Younghae Do and Jaeun Lee for their attention to this work.
The \TeX\, macro package PSTricks~\cite{PSTricks} was essential for
typesetting the equations and figures.

\bibliographystyle{amsplain}

\end{document}